\documentclass[11pt]{article}
\usepackage{amsmath}
\usepackage{amssymb,amsfonts,amsmath,latexsym,cite,color, psfrag}
\newtheorem{theo}{Theorem}[section]

\newtheorem{lem} [theo]{Lemma}
\newtheorem{coro}[theo]{Corollary}
\newtheorem{prop}[theo]{Proposition}

\newtheorem{conjecture}[theo]{Conjecture}
\def\R{{\widetilde{R}}}
\allowdisplaybreaks

\makeatletter \@addtoreset{equation}{section}
\@addtoreset{theo}{}\makeatother

\def\qed{\hfill \rule{4pt}{7pt}}
\def\pf{\noindent {\it Proof.} }

\def\R{{\widetilde{R}}}
\begin{document}
\begin{center}
{\Large\bf A Class of Kazhdan-Lusztig $R$-Polynomials and

 $q$-Fibonacci Numbers}
\end{center}
 \vskip 4mm
 \begin{center}
{\small William Y.C. Chen$^1$, Neil J.Y. Fan$^2$, Peter L. Guo$^3$, Michael X.X. Zhong$^4$}

\vskip 4mm
$^{1,3,4}$Center for Combinatorics, LPMC-TJKLC\\
Nankai University,
Tianjin 300071,
P.R. China \\[3mm]

$^2$Department of Mathematics\\
Sichuan University, Chengdu, Sichuan 610064, P.R. China

\vskip 4mm

$^1$chen@nankai.edu.cn,  $^2$fan@scu.edu.cn,
$^3$lguo@nankai.edu.cn\\
$^4$michaelzhong@mail.nankai.edu.cn
\end{center}

\begin{abstract}
Let $S_n$ denote the symmetric group on $\{1,2,\ldots,n\}$.
For two permutations $u, v\in S_n$ such that $u\leq v$ in the Bruhat order, let
$R_{u,v}(q)$ and  $\R_{u,v}(q)$ denote the Kazhdan-Lusztig $R$-polynomial and
$\R$-polynomial, respectively.
Let $v_n=34\cdots n\, 12$, and let $\sigma$ be a permutation
such that $\sigma\leq v_n$.
We obtain a formula for the
$\R$-polynomials $\R_{\sigma,v_n}(q)$ in terms of the $q$-Fibonacci numbers
depending on a parameter determined by the
reduced expression of $\sigma$.
When $\sigma$ is the identity $e$, this reduces to a formula obtained by Pagliacci.
In another direction, we obtain a formula for the $\R$-polynomial
$\R_{e,\,v_{n,i}}(q)$, where  $v_{n,i} = 3 4\cdots i\,n\, (i+1)\cdots (n-1)\, 12$.
In a more general context, we conjecture that for any two permutations $\sigma,\tau\in S_n$ such that
 $\sigma\leq \tau\leq v_n$,   the $\R$-polynomial
$\R_{\sigma,\tau}(q)$ can be expressed as a product of   $q$-Fibonacci numbers
multiplied by a power of $q$.
\end{abstract}



\section{Introduction}

Let $S_n$ denote the symmetric group on $\{1,2,\ldots,n\}$.
For two permutations $u, v\in S_n$ such that $u\leq v$ in the Bruhat order, let
$R_{u,v}(q)$ be the Kazhdan-Lusztig $R$-polynomial, and $\R_{u,v}(q)$ be
the  Kazhdan-Lusztig $\R$-polynomial.
Let $v_n=34\cdots n\, 12$, and let $\sigma$ be a permutation
such that $\sigma\leq v_n$.
The main result of this paper is a formula for the
$\R$-polynomials $\R_{\sigma,\,v_n}(q)$ in terms of the $q$-Fibonacci numbers
depending on a parameter determined by the
reduced expression of $\sigma$.
When $\sigma$ is the identity
permutation $e$, a formula for the $\R$-polynomials has been given by
Pagliacci \cite[Theorem 4.1]{Pagliacci1}.

We also derive a formula for the $\R$-polynomials
$\R_{e,\,v_{n,i}}(q)$, where  $v_{n,i} = 3 4\cdots i\,n\, (i+1)\cdots (n-1)\, 12$, which can be viewed as
a generalization of Pagliacci's formula \cite[Theorem 4.1]{Pagliacci1} in another direction.
We conclude this paper with a  conjecture that for any two permutations $\sigma,\tau\in S_n$ such that
 $\sigma\leq \tau\leq v_n$,   the $\R$-polynomial
$\R_{\sigma,\,\tau}(q)$ can be expressed as a product of $q$-Fibonacci numbers
and a power of $q$.

Let us give an overview of some notation and background.
For each permutation $\pi$ in $S_n$, it is known that $\pi$ can be
expressed as a product of simple transpositions $s_i=(i,i+1)$
subject to the following braid relations
\begin{align*}
s_i\,s_j&=s_j\,s_i,\ \ \ \ \ \ \ \ \ \  \ \ \ \ \ \ {\rm for}\ |i-j|>1;\\[5pt]
s_i\,s_{i+1}\,s_i&=s_{i+1}\,s_i\,s_{i+1},\ \ \ \ \ \ \  \ {\rm for}\ 1 \leq i\leq n-2.
\end{align*}
An expression $\omega$ of $\pi$ is said to be reduced if the number of
simple transpositions appearing  in $\omega$  is minimum.
The following word property is due to Tits, see Bj\"orner and Brenti \cite[Theorem 3.3.1]{BB}.

\begin{theo}[\mdseries{Word Property}]\label{Tits}
Let $\pi$ be a permutation of $S_n$, and  $\omega_1$ and $\omega_2$ be two reduced expressions of $\pi$.
Then $\omega_1$ and $\omega_2$ can be obtained from each other by applying a sequence of braid relations.
\end{theo}

Let $\ell(\pi)$ denote the length of $\pi$, that is, the number of simple transpositions in a reduced expression of $\pi$.
Write  $D_R(\pi)$ for the set of right descents of $\pi$, namely,
\[
D_R(\pi)=\{s_i\colon 1\leq i\leq n-1,\, \ell(\pi s_i)<\ell(\pi)\}.
\]

The exchange condition
gives a characterization for the  (right) descents of a permutation
in terms of reduced expressions, see Humphreys \cite[Section 1.7]{Humphreys}.

\begin{theo}[\mdseries{Exchange Condition}]\label{exchange}
Let $\pi=s_{i_1}s_{i_2}\cdots s_{i_k}$ be a reduced expression of $\pi$. If $s_i\in D_R(\pi)$,
then there exists an
index $i_j$ for which  $\pi s_i=s_{i_1}\cdots \widehat{s}_{i_j}\cdots s_{i_k}$, where $\widehat{s}_{i_j}$
means that $s_{i_j}$ is missing.
In particular, $\pi$ has a reduced expression  ending with $s_i$ if and only if $s_i\in D_R(\pi)$.
\end{theo}

The following subword property serves as a definition of the Bruhat order.
For other equivalent definitions of the Bruhat order, see Bj\"orner and Brenti
\cite{BB}.
For a reduced expression  $\omega=s_{i_1}s_{i_2}\cdots s_{i_k}$, we say that $s_{i_{j_1}}s_{i_{j_2}}\cdots s_{i_{j_m}}$ is a
subword of $\omega$ if $1\leq j_1<j_2<\cdots < j_m\leq k$.

\begin{theo}[\mdseries{Subword Property}]\label{subword}
Let $u$ and $v$ be two permutation in $S_n$.
Then $u\leq v$ in the Bruhat order if and only if every reduced expression of $v$
has a subword that is  a reduced expression  of $u$.
\end{theo}

The Bruhat order satisfies  the following lifting property, see Bj\"orner and Brenti
 \cite[Proposition 2.2.7]{BB}.

\begin{theo}[\mdseries{Lifting Property}]\label{lift}
Suppose that $u$ and $v$ are two permutations in $S_n$ such that $u<v$. For any
 simple transposition
 $s_i$ in $D_R(v)\backslash D_R(u)$,
 we have $u\leq v s_i$ and $us_i\leq v$.
\end{theo}

The Kazhdan-Lusztig $R$-polynomials, which
were introduced by Kazhdan and Lusztig \cite{KL1},
can be recursively determined  by the following properties, see also Humpreys \cite[Section 7.5]{Humphreys}.

\begin{theo}
For any $u,v\in S_n$,
\begin{itemize}
\item[$\mathrm{(i)}$] $R_{u,v}(q)=0$, if $u\nleq v$;
\item[$\mathrm{(ii)}$] $R_{u,v}(q)=1$, if $u= v$;
\item[$\mathrm{(iii)}$] If $u<v$ and $s\in D_{R}(v)$,
\[R_{u,v}(q)=\left\{
        \begin{array}{ll}
          R_{us,vs}(q), & \hbox{\rm{if} $s\in D_{R}(u)$;} \\[5pt]
          qR_{us,vs}(q)+(q-1)R_{u,vs}(q), & \hbox{\rm{if} $s \notin  D_{R}(u)$.}
        \end{array}
      \right.
\]
\end{itemize}
\end{theo}

While $R$-polynomials may contain   negative coefficients,
a variant of the $R$-polynomials introduced by Dyer \cite{Dyer}, which
 has been called the $\widetilde{R}$-polynomials,
has only nonnegative coefficients. For an alternative definition of the
$\R$-polynomials for the symmetric group, see Brenti \cite{Brenti2}.
The following two theorems are due to Dyer \cite{Dyer}, see also
Brenti \cite{Brenti2}.

\begin{theo}\label{Rtilde}
Let   $u,v\in S_n$ with $u\leq v$. Then, for  $s\in D_R(v)$,
\begin{equation}\label{def1}
\widetilde{R}_{u,v}(q)=\left\{
    \begin{array}{ll}
      \widetilde{R}_{us,vs}(q), & \hbox{\rm{if} $s\in D_R(u)$;} \\[5pt]
      \widetilde{R}_{us,vs}(q)+q\widetilde{R}_{u,vs}(q), & \hbox{\rm{if} $s\notin D_R(u)$.}
    \end{array}
  \right.
\end{equation}
\end{theo}

\begin{theo}\label{relation}
Let   $u,v\in S_n$ with $u\leq v$. Then
\[
R_{u,v}(q)
=q^{\frac{\ell(v)-\ell(u)}{2}}\widetilde{R}_{u,v}(q^{\frac{1}{2}}-q^{-\frac{1}{2}}).
\]
\end{theo}

Recall that $v_n=34\cdots n\, 12$ and $v_{n,i} = 3 4\cdots i\,n\, (i+1)\cdots (n-1)\, 12$. We shall use the recurrence relations in Theorem \ref{Rtilde} to  deduce  a formula for
 the $\R$-polynomials $\R_{\sigma,\,v_n}(q)$, from which we
 also find a formula for the $\R$-polynomial
$\R_{e,\,v_{n,i}}(q)$.

\section{Main result}

 The main result of this paper is an equation for $\R_{\sigma,\, v_n}(q)$,
 where $v_n=34\cdots n\,12$ and $\sigma\leq v_n$ in the Bruhat order.
Combining this equation with a formula of Pagliacci \cite{Pagliacci1}, we obtain
an expression of $\R_{\sigma,\,v_n}(q)$ in terms of $q$-Fibonacci numbers.
 To describe our result, we need the following reduced expression of $v_n$.

\begin{prop}
For
$n\geq 3$,
\begin{equation*}\label{nf}
\Omega_n=s_2s_1s_3s_2\,\cdots\,s_{n-1}s_{n-2}
\end{equation*}
is a reduced expression of $v_n$.
\end{prop}

Let  $\sigma$ be a permutation of $S_n$ such that $\sigma\leq v_n$.
By the subword property in Theorem \ref{subword},
  $\sigma$ can be expressed as a reduced subword of $\Omega_n$.
We  introduce two
statistics of a reduced subword of $\Omega_n$.

Let $\omega=s_{i_1}s_{i_2}\cdots s_{i_k}$ be a reduced  subword  of $\Omega_n$.
Define
\begin{align*}
D(\omega)&=\{1\leq t<k \colon i_t-i_{t+1}= 1\}.
\end{align*}
We use $d(\omega)$ to denote the cardinality of $D(\omega)$, and let
\begin{equation}\label{relation-x}
h(\omega)=n-\ell(\omega)+d(\omega).
\end{equation}
For example, for a reduced subword $\omega=s_2s_3s_4s_3 s_6s_5s_7$ of $\Omega_9$, we have
$D(\omega)=\{3,5\}$, and thus
$d(\omega)=2$ and $h(\omega)=4$.
Note that  $h(\omega)$ depends on both $\omega$  and  $n$.
This causes no confusion since the index $n$
is always clear from the context.

The main result in this paper is the following equation  for the $\R$-polynomials $\R_{\sigma,\,v_n}(q)$.

\begin{theo}\label{main}
For $n\geq 3$, let $\sigma$ be a permutation in $S_n$ such that $\sigma\leq v_n$, and let $\omega$
be any  reduced expression of $\sigma$ that is a subword of $\Omega_n$.
 Then we have
\begin{equation}\label{findittoday-33}
\R_{\sigma,\,v_n}(q)=q^{\ell(\omega)-2d(\omega)}\R_{e,\, v_{h(\omega)}}(q).
\end{equation}
\end{theo}

Let $F_n(q)$ be the $q$-Fibonacci numbers, that is,   $F_0(q)=F_1(q)=1$
and for $n\geq 2$,
\[F_n(q)=F_{n-1}(q)+qF_{n-2}(q).\]
Pagliacci \cite[Theorem 4.1]{Pagliacci1} has shown that
\begin{equation}\label{zxbn}
\R_{e,\, v_n}(q)=q^{2n-4}F_{n-2}(q^{-2}).
\end{equation}
As a consequence of Theorem \ref{main}  and formula \eqref{zxbn}, we
obtain an expression of $\R_{\sigma,\,v_n}(q)$ in terms of $q$-Fibonacci numbers.

\begin{coro} \label{coro-1}
For $n\geq 3$, let $\sigma$ be a permutation in $S_n$ such that $\sigma\leq v_n$, and let $\omega$
be any  reduced expression of $\sigma$ that is a subword of $\Omega_n$.
 Then we have
\begin{equation}\label{findittoday-1}
\R_{\sigma,\,v_n}(q)=q^{2n-\ell(\sigma)-4} F_{h(\omega)-2}(q^{-2}).
\end{equation}
\end{coro}

To give an inductive proof of  Theorem \ref{main}, we need three lemmas.
Assume that $\omega$ is a reduced subword  of $\Omega_n$.
When $\omega s_{n-1} \leq \Omega_n$, the first two lemmas are
 concerned with the existence of a reduced expression $\omega'$ of $\omega s_{n-1}$ such that $d(\omega')=d(\omega)$.
When  $\omega s_{n-1} \not\leq \Omega_n$, the third lemma shows  that $h(\omega)=2$.

\begin{lem}\label{lemma-1}
Let $\omega$ be a reduced subword of $\Omega_{n}$. If $\omega s_{n-1}\leq \Omega_n$ and  $s_{n-1}\in D_R(\omega)$, then
there exists a reduced expression $\omega'$ of $\omega s_{n-1}$ such that
$\omega'$ is a subword of $\Omega_{n}$ and
$d(\omega')=d(\omega)$.
\end{lem}

\pf  We use induction on $n$.
It is easy to check that the lemma holds for $n\leq 3$. Assume
that $n>3$ and the assertion holds for $n-1$.
We now consider the case for $n$.
By definition, we have $\Omega_n=\Omega_{n-1}s_{n-1}s_{n-2}$.
Since $\omega$ is a subword of $\Omega_{n}$, we can write
$\omega=\omega_1 \omega_2$, where $\omega_1$
is a subword  of $\Omega_{n-1}$ and  $\omega_2$ is a subword of $s_{n-1}s_{n-2}$.
 Because   $s_{n-1}\in D_R(\omega)$,
we have the following
 two cases.

\noindent
Case 1: $\omega=\omega_1 s_{n-1}$. Set $\omega'=\omega_1$.
Clearly,  $\omega'$   is a reduced expression of $\omega s_{n-1}$.
Moreover, it is easy to check that $D(\omega)=D(\omega')$, and thus $d(\omega')=d(\omega)$, that is,   $\omega'$ is a desired reduced expression of $\omega s_{n-1}$.

\noindent
Case 2: $\omega=\omega_1 s_{n-1}s_{n-2}$.
Since $s_{n-1}\in D_R(\omega)$,
by Theorem \ref{exchange}, there exists a  reduced expression
of $\omega$ ending with $s_{n-1}$.
Hence the word property in  Theorem \ref{Tits} ensures
  the existence of a reduced expression of $\omega_1$ ending with $s_{n-2}$.
This implies that $s_{n-2}$ belongs to $D_R(\omega_1)$.
By the induction hypothesis, there exists a reduced expression $\omega_1'$
of $\omega_1s_{n-2}$ such that $\omega_1'$ is a subword of $\Omega_{n-1}$ and
$d(\omega_1')=d(\omega_1)$.

Set
\[\omega'=\omega_1's_{n-1}s_{n-2}.\]
We deduce that $\omega'$ is a desired reduced subword.
Since
\begin{equation*}\label{aj}
\omega'=\omega_1' s_{n-1}s_{n-2}=\omega_1 s_{n-2}s_{n-1}s_{n-2}=\omega_1 s_{n-1}s_{n-2}s_{n-1}=\omega s_{n-1},
\end{equation*}
we see that $\omega'$ is an expression of $\omega s_{n-1}$.
On the other hand, since $\omega'$ consists of $\ell(\omega_1')+2$
simple transpositions and
\[\ell(\omega_1')+2=\ell(\omega_1)+1=\ell(\omega)-1=
\ell(\omega s_{n-1}),\]
we conclude that $\omega'$ is a reduced  expression of $\omega s_{n-1}$.
By the construction of $\omega'$,
we have
\[d(\omega')=d(\omega_1')+1=d(\omega_1)+1=d(\omega).\]
This completes the proof. \qed

The next lemma deals with the case $s_{n-1} \not\in D_R(\omega)$.

\begin{lem}\label{lemma-3}
Let $\omega$ be a reduced  subword of $\Omega_{n}$. If $\omega s_{n-1}\leq \Omega_n$ and $s_{n-1}\not\in D_R(\omega)$, then
there exists a reduced expression $\omega'$ of $\omega s_{n-1}$ such that
$\omega'$ is a subword of $\Omega_{n}$  and
$d(\omega')=d(\omega)$.

\end{lem}

\pf We use induction on $n$.
It is easily checked that the lemma holds for $n\leq 3$. Assume
that $n>3$ and the assertion holds for $n-1$.
We now consider the case for $n$.
Let
$\omega=\omega_1 \omega_2$, where $\omega_1$
is a subword  of $\Omega_{n-1}$ and  $\omega_2$ is a subword of $s_{n-1}s_{n-2}$.
Since   $s_{n-1}\not\in D_R(\omega)$,
we have the following three cases.

\noindent
Case 1: $\omega=\omega_1$.
Set  $\omega'=\omega_1 s_{n-1}$.  It is easily seen that $\omega'$ is a desired reduced expression.

\noindent
Case 2: $\omega=\omega_1 s_{n-2}$.
 We claim that $\omega=\omega_1 s_{n-2}\leq \Omega_{n-1}$.
Note that $s_{n-1}$ does not appear in $\omega_1$. Since $\omega s_{n-1}=\omega_1 s_{n-2}s_{n-1}$, by   Theorem \ref{Tits}, there does not exist any reduced expression of
$\omega s_{n-1}$ ending with
 $s_{n-2}$. This implies that $s_{n-2}$ does not belong to $D_R(\omega s_{n-1})$.
Thus, by the lifting property in Theorem \ref{lift},
we deduce that
\[\omega s_{n-1}\leq \Omega_n s_{n-2}=\Omega_{n-1} s_{n-1}.\]
This implies that  $\omega=\omega_1  s_{n-2}\leq \Omega_{n-1}$, as claimed.

Since $\omega=\omega_1  s_{n-2}$ is reduced, we see that $s_{n-2}\not\in D_R(\omega_1)$. By the induction hypothesis,
there exists a reduced expression $\omega_1'$ of $\omega_1  s_{n-2}$ such that
$\omega_1'$ is a subword of $\Omega_{n-1}$  and
$d(\omega_1')=d(\omega_1)$.
Set $\omega'=\omega_1'  s_{n-1}$. We find that
$\omega'$ is a   reduced expression of $\omega  s_{n-1}$
such that $d(\omega')=d(\omega)$.

\noindent
Case 3: $\omega=\omega_1  s_{n-1}s_{n-2}$.
We claim  that $s_{n-2}\not\in D_R(\omega_1)$.
Suppose to the contrary that $s_{n-2}\in D_R(\omega_1)$.
By Theorem \ref{exchange}, there exists a reduced expression  of $\omega_1$
ending with $s_{n-2}$.  Write $\omega_1=\mu s_{n-2}$, where $\mu$ is a reduced expression.
Then we get
\[\omega=\mu s_{n-2}s_{n-1}s_{n-2}=\mu s_{n-1}s_{n-2}s_{n-1},\]
contradicting  the assumption that $s_{n-1}\not\in D_R(\omega)$.
So the claim is proved.

On the other hand, since
\[\omega_1  s_{n-2}s_{n-1}s_{n-2}=\omega_1  s_{n-1}s_{n-2}s_{n-1}=\omega s_{n-1}\leq \Omega_n,\]
we have $\omega_1 s_{n-2}\leq \Omega_{n-1}$.
It follows from the induction hypothesis that
there exists a reduced expression $\omega_1'$ of $\omega_1  s_{n-2}$ such that
$\omega_1'$ is a subword of $\Omega_{n-1}$  and
$d(\omega_1')=d(\omega_1)$.

Let \[\omega'=\omega_1'  s_{n-1}s_{n-2}.\]
Since
\[\omega'=\omega_1'  s_{n-1}s_{n-2}=\omega_1  s_{n-2}s_{n-1}s_{n-2}=\omega_1  s_{n-1}s_{n-2}s_{n-1}=\omega  s_{n-1},\]
we deduce that $\omega'$ is a  reduced expression  of $\omega  s_{n-1}$.
By the construction of
$\omega'$,
we obtain that
\[d(\omega')=d(\omega_1')+1=d(\omega_1)+1=d(\omega),\]
as required.
\qed

We now come to the third lemma.

\begin{lem}\label{lemma-2}
Let $\omega$ be a reduced  subword of $\Omega_{n}$.
If $\omega  s_{n-1}\not\leq \Omega_{n}$, then we have $h(\omega)=2$.
\end{lem}

\pf
We proceed by induction on $n$.
It can be verified that the lemma holds for $n\leq 3$. Assume
that $n>3$ and the assertion holds for $n-1$.
Consider the case for $n$.
Write
$\omega=\omega_1  \omega_2$, where $\omega_1$
is a subword  of $\Omega_{n-1}$ and  $\omega_2$ is a subword of $s_{n-1}s_{n-2}$.
Since $\omega  s_{n-1}\not\leq \Omega_{n}$, we see that
$s_{n-1}$ is not a right descent of $\omega$. We have the following two cases.

\noindent
Case 1: $\omega=\omega_1 s_{n-2}$.
Since $\omega  s_{n-1}\not\leq \Omega_{n}$, we have
$\omega=\omega_1  s_{n-2}\not\leq \Omega_{n-1}$. Thus, by the induction hypothesis,
we get $h(\omega_1)=2$. Noticing that $\ell(\omega_1)=\ell(\omega)-1$ and
$d(\omega_1)=d(\omega)$, we obtain that
\[h(\omega)=n-\ell(\omega)+d(\omega)=
n-1-\ell(\omega_1)+d(\omega_1)=h(\omega_1)=2,\]
as required.

\noindent
Case 2: $\omega=\omega_1 s_{n-1}s_{n-2}$. We  claim that $\omega_1  s_{n-2}\not\leq \Omega_{n-1}$. Suppose to the contrary that $\omega_1  s_{n-2}\leq \Omega_{n-1}$.
Note that
\[\omega s_{n-1}=\omega_1 s_{n-1}s_{n-2}s_{n-1}
=\omega_1 s_{n-2}s_{n-1}s_{n-2}.\]
This yields $\omega s_{n-1}\leq \Omega_n$, contradicting the assumption
that $\omega s_{n-1}\not\leq \Omega_n$. So the claim is verified.

By the induction hypothesis, we have $h(\omega_1)=2$.
Since $\ell(\omega_1)=\ell(\omega)-2$ and
$d(\omega_1)=d(\omega)-1$,  we find that
\[h(\omega)=n-\ell(\omega)+d(\omega)
=n-1-\ell(\omega_1)+d(\omega_1)=h(\omega_1)=2,\]
as requied.\qed

Next we give a proof of  Theorem \ref{main} based on the
above lemmas.

\noindent\textit{Proof of Theorem \ref{main}. }
Let
\[T_n(q)=\R_{e,\,v_n}(q).\]
and
\begin{equation}\label{relation-2}
g(\omega)=\ell(\omega)-2d(\omega).
\end{equation}
Then equation \eqref{findittoday-33}   can be rewritten as
\begin{equation}\label{findittoday-3}
\R_{\sigma,\,v_n}(q)=q^{g(\omega)}T_{h(\omega)}(q).
\end{equation}

We  proceed to prove  \eqref{findittoday-3}   by induction on $n$.
It can be checked that
\eqref{findittoday-3} holds for $n\leq 3$.
Assume that $n>3$ and \eqref{findittoday-3} holds for $n-1$.
For the case for $n$, let $\omega=\omega_1  \omega_2$, where $\omega_1$
is a subword  of $\Omega_{n-1}$ and  $\omega_2$ is a subword of $s_{n-1}s_{n-2}$.
There are four cases.

\noindent
Case 1: $\omega=\omega_1 s_{n-2}$.  It follows from \eqref{def1}  that
\begin{align}
\R_{\sigma,\,v_n}(q)&=\R_{\omega_1 s_{n-2},\,\Omega_{n-1}s_{n-1}s_{n-2}}(q)\nonumber\\[5pt]
&=\R_{\omega_1,\,\Omega_{n-1}s_{n-1}}(q)\nonumber\\[5pt]
&=\R_{\omega_1  s_{n-1},\,\Omega_{n-1}}(q)+q\R_{\omega_1,\,\Omega_{n-1}}(q).\label{sgh-1}
\end{align}
Since $\omega_1  s_{n-1}\not\leq\Omega_{n-1}$, we see that
the first term in \eqref{sgh-1} vanishes. Thus, \eqref{sgh-1} becomes
\begin{equation}\label{sgh-11}
\R_{\sigma,\,v_n}(q)=q\R_{\omega_1,\,\Omega_{n-1}}(q).
\end{equation}
By the induction hypothesis, we have
\begin{equation}\label{sgh-2}
\R_{\omega_1,\,\Omega_{n-1}}(q)=q^{g(\omega_1)}T_{h(\omega_1)}.
\end{equation}
Since $d(\omega_1)=d(\omega)$,  we get
\begin{equation}\label{t1}
g(\omega_1)=\ell(\omega_1)-2d(\omega_1)=\ell(\omega)-1-2d(\omega)=g(\omega)-1
\end{equation}
and
\begin{equation}\label{t2}
h(\omega_1)=n-1-\ell(\omega_1)+d(\omega_1)=n-\ell(\omega)+d(\omega)=h(\omega).
\end{equation}
Plugging \eqref{t1} and \eqref{t2} into \eqref{sgh-2}, we obtain
\[\R_{\omega_1,\,\Omega_{n-1}}(q)=q^{g(\omega_1)}T_{h(\omega_1)}
=q^{g(\omega)-1}T_{h(\omega)}(q),\]
which leads to
\[\R_{\sigma,\,v_n}(q)=q\R_{\omega_1,\,\Omega_{n-1}}(q)
=q^{g(\omega)}T_{h(\omega)}(q).\]

\noindent
Case 2: $\omega=\omega_1 s_{n-1}$. By Theorem \ref{Tits}, there is no
 reduced expression of
$\sigma$ that ends with $s_{n-2}$. This implies that $s_{n-2}$
is not a right descent of $\sigma$, so that
\begin{align}
\R_{\sigma,\,v_n}(q)&=\R_{\omega_1 s_{n-1}s_{n-2},\,\Omega_{n-1} s_{n-1}}(q)+
q\R_{\omega_1 s_{n-1},\,\Omega_{n-1}s_{n-1}}(q)\nonumber\\[5pt]
&=\R_{\omega_1 s_{n-1}s_{n-2},\,\Omega_{n-1} s_{n-1}}(q)+q\R_{\omega_1,\,\Omega_{n-1}}(q).\label{formula-1}
\end{align}

We claim that the first term in \eqref{formula-1} vanishes, or equivalently,
 $\omega_1 s_{n-1}s_{n-2}\not\leq \Omega_{n-1}s_{n-1}$.
 Suppose to the contrary that $\omega_1 s_{n-1}s_{n-2}\leq \Omega_{n-1}s_{n-1}$.
 By Theorem \ref{subword}, there exists a subword $\mu$ of $\Omega_{n-1}s_{n-1}$
that is a reduced expression of $\omega_1 s_{n-1}s_{n-2}$. Since
$s_{n-1}$ must appear in $\mu$, we may write  $\mu$ in the following form
\begin{equation*}\label{relation-12}
\mu=s_{i_1}s_{i_2}\cdots s_{i_{k}}  s_{n-1},
\end{equation*}
 where $s_{i_1}s_{i_2}\cdots s_{i_{k}}$ is a reduced subword of $\Omega_{n-1}$.
By the word property in Theorem \ref{Tits},   $\omega_1 s_{n-1}s_{n-2}$ can be
obtained from  $\mu$ by applying the braid relations.
However, this is impossible since any simple transposition $s_{n-2}$ appearing in $\mu$
cannot be moved to the last position by applying the braid relations.
So the claim is proved, and hence  \eqref{formula-1} becomes
\[
\R_{\sigma,\,v_n}(q)
=q\R_{\omega_1,\,\Omega_{n-1}}(q).\]
It is easily seen that
 \[g(\omega_1)=g(\omega)-1\ \ \text{and}\ \ h(\omega_1)=h(\omega).\]
By the
induction hypothesis, we deduce that
\[
\R_{\sigma,\,v_n}(q)
=q\R_{\omega_1,\,\Omega_{n-1}}(q)
= q^{g(\omega_1)+1}T_{h(\omega_1)}(q)
=q^{g(\omega)}T_{h(\omega)}(q).
\]

\noindent
Case 3: $\omega=\omega_1 s_{n-1}s_{n-2}$. It is clear from \eqref{def1} that
\begin{equation}\label{aohh}
\R_{\sigma,\,v_n}(q)=\R_{\omega_1s_{n-1},\,\Omega_{n-1}s_{n-1}}(q)=\R_{\omega_1,\,\Omega_{n-1}}(q).
\end{equation}
Noting that $d(\omega)=d(\omega_1)+1$, we obtain
\[g(\omega_1)=\ell(\omega_1)-2d(\omega_1)=\ell(\omega_1)+2-2d(\omega)=\ell(\omega)-2d(\omega)=g(\omega)\]
and
\[h(\omega_1)=n-1-\ell(\omega_1)+d(\omega_1)=n-2-\ell(\omega_1)+d(\omega)=n-\ell(\omega)+d(\omega)=h(\omega).\]
Thus, by \eqref{aohh} and  the induction hypothesis, we find that
\[
\R_{\sigma,\,v_n}(q)=\R_{\omega_1,\,\Omega_{n-1}}(q)
=q^{g(\omega_1)}T_{h(\omega_1)}(q)=q^{g(\omega)}T_{h(\omega)}(q).
\]

\noindent
Case 4: $\omega=\omega_1$. Here are two subcases.

\noindent
Subcase 1: $s_{n-2}\in D_R(\omega_1)$. By \eqref{def1}, we deduce that
\begin{align}
\R_{\sigma,\, v_n}(q)&=\R_{\omega_1 s_{n-2},\, \Omega_{n-1}s_{n-1}}(q)\nonumber \\
&=q\R_{\omega_1 s_{n-2},\,\Omega_{n-1}}(q) .\label{ckj}
\end{align}
By Lemma \ref{lemma-1}, there exists a reduced expression $\omega_1'$ of $\omega_1 s_{n-2}$ such that
$\omega_1'$ is a subword of $\Omega_{n-1}$ and $d(\omega_1')=d(\omega_1)$.
Consequently,
\begin{equation*}
g(\omega_1')=\ell(\omega_1')-2d(\omega_1')=
\ell(\omega_1)-1-2d(\omega_1)=g(\omega_1)-1
\end{equation*}
and
\begin{equation*}
h(\omega_1')=n-1-\ell(\omega_1')+d(\omega_1')=
n-\ell(\omega_1)+d(\omega_1)=h(\omega_1)+1.
\end{equation*}
By  the induction hypothesis, we obtain that
\begin{align}
\R_{\omega_1 s_{n-2},\,\Omega_{n-1}}(q)&=\R_{\omega_1' ,\, \Omega_{n-1}}(q)\nonumber \\[5pt]
&=q^{g(\omega_1')}T_{h(\omega_1')}(q)\nonumber\\[5pt]
&=q^{g(\omega_1)-1}T_{h(\omega_1)+1}(q).\label{aln-2}
\end{align}
But  $h(\omega)=h(\omega_1)+1$,
substituting \eqref{aln-2} into  \eqref{ckj} gives
\begin{align*}
\R_{\sigma,\,v_n}(q)
&=q\R_{\omega_1 s_{n-2},\,\Omega_{n-1}}(q)\\[5pt]
&=q^{g(\omega)}T_{h(\omega)}(q).
\end{align*}

\noindent
Subcase 2: $s_{n-2}\not\in D_R(\omega_1)$. By \eqref{def1}, we see that
\begin{align}
\R_{\sigma,\,v_n}(q)&=\R_{\omega_1  s_{n-2},\,\Omega_{n-1}s_{n-1}}(q)+q\R_{\omega_1 ,\,\Omega_{n-1}s_{n-1}}(q)\nonumber\\[5pt]
&=q\R_{\omega_1  s_{n-2},\,\Omega_{n-1}}(q)+q^2\R_{\omega_1,\,\Omega_{n-1}}(q)\label{findtitoday-2}
\end{align}
By the induction hypothesis, the second term in \eqref{findtitoday-2} equals
\begin{equation}\label{nnn-1}
q^2\R_{\omega_1 ,\,\Omega_{n-1}}(q)=q^{g(\omega_1)+2}T_{h(\omega_1)}(q)=q^{g(\omega)+2}T_{h(\omega)-1}(q).
\end{equation}
It remains to compute  the first term in \eqref{findtitoday-2}. To this end,
we have the following two cases.

\noindent
Subcase 2a: $\omega_1  s_{n-2}\leq \Omega_{n-1}$.
By Lemma \ref{lemma-3},
there exists a reduced expression $\omega_1'$ of $\omega_1  s_{n-2}$ such that
$\omega_1'$ is a subword of $\Omega_{n-1}$ and
$d(\omega_1)=d(\omega_1')$.
Hence
\begin{equation*}
g(\omega_1')=\ell(\omega_1')-2d(\omega_1')=
\ell(\omega_1)+1-2d(\omega_1)=g(\omega_1)+1=
g(\omega)+1
\end{equation*}
and
\begin{equation*}
h(\omega_1')=n-1-\ell(\omega_1')+d(\omega_1')=
n-2-\ell(\omega_1)+d(\omega_1)=h(\omega_1)-1=
h(\omega)-2.
\end{equation*}
By the induction hypothesis, we obtain that
\begin{equation}\label{nal}
\R_{\omega_1  s_{n-2},\,\Omega_{n-1}}(q)= q^{g(\omega_1')}T_{h(\omega_1')}(q)
=q^{g(\omega)+1}T_{h(\omega)-2}(q).
\end{equation}

Putting \eqref{nnn-1} and \eqref{nal} into \eqref{findtitoday-2},
we deduce that
\begin{align}\label{re-30}
\R_{\sigma,v_n}(q)
&=q\R_{\omega_1 s_{n-2},\,\Omega_{n-1}}(q)+q^2\R_{\omega_1,\,\Omega_{n-1}}(q) \nonumber \\[5pt]
&=q^{g(\omega)+2}T_{h(\omega)-2}(q)+q^{g(\omega)+2}T_{h(\omega)-1}(q)\nonumber\\[5pt]
&=q^{g(\omega)}(q^2T_{h(\omega)-2}(q)+q^2T_{h(\omega)-1}(q)).
\end{align}
In view of the following relation due to Pagliacci \cite{Pagliacci1}
\[T_n(q)=q^2T_{n-2}(q)+q^2T_{n-1}(q),\]
 \eqref{re-30} can be rewritten as
\[\R_{\sigma,\,v_n}(q)=q^{g(\omega)}T_{h(\omega)}(q).\]

\noindent
Subcase 2b: $\omega_1 s_{n-2}\not\leq \Omega_{n-1}$.
In this case, we have
\begin{equation}\label{p1}
 \R_{\omega_1 s_{n-2},\,\Omega_{n-1}}(q)=0.
 \end{equation}
By Lemma \ref{lemma-2}, we find that $h(\omega_1)=2$.
Thus \eqref{nnn-1} reduces to
 \begin{equation}\label{p2}
 q^2\R_{\omega_1 ,\,\Omega_{n-1}}(q)=q^{g(\omega)+2}T_2(q)=q^{g(\omega)+2}.
 \end{equation}
Putting \eqref{p1} and  \eqref{p2}  into \eqref{findtitoday-2}, we get
\begin{equation}\label{hoag}
\R_{\sigma,\,v_n}(q)=q^{g(\omega)+2}.
\end{equation}
Since $T_3(q)=q^2$ and $h(\omega)=h(\omega_1)+1=3$, it follows from \eqref{hoag}  that
\[\R_{\sigma,\,v_n}(q)=q^{g(\omega)}T_3(q)=q^{g(\omega)}T_{h(\omega)}(q),\]
and hence the proof is complete. \qed

 For $2\leq i\leq n-1$, let
 \[v_{n,i}=\left\{
        \begin{array}{ll}
          n3 4\cdots(n-1)\, 12, & \hbox{\rm{if} $i=2$;} \\[5pt]
          3 4\cdots i\,n\, (i+1)\cdots (n-1)\, 12, & \hbox{\rm{if} $3\leq i\leq n-1$.}
        \end{array}
      \right.
\]
We obtain the following formula for
 $\R_{e,\, v_{n,i}}(q)$, which reduces to formula
\eqref{zxbn} due to Pagliacci in the case $i=n-1$.

\begin{theo}\label{main2}
Let  $n\geq 3$ and $2\leq i\leq n-1$. Then we have
\begin{equation}\label{0930-1}
\R_{e,\, v_{n,i}}(q)=\sum_{k=0}^{n-i-1}q^{3n-i-2k-5}{n-i-1 \choose k}F_{n-k-2}(q^{-2}).
\end{equation}
\end{theo}

\pf Recall that
\[
T_n(q)=\R_{e,\, v_n}(q)=q^{2n-4}F_{n-2}(q^{-2}),
\]
see \eqref{zxbn}.
Hence   \eqref{0930-1} can be rewritten as
\begin{equation}\label{0930-t}
\R_{e,\, v_{n,i}}(q)=\sum_{k=0}^{n-i-1}q^{n-i-1}{n-i-1 \choose k}T_{n-k}(q).
\end{equation}
Note that $\Omega_n s_{n-3}s_{n-4}\cdots s_{i-1}$
is a reduced expression of the permutation $v_{n,i}$. By the defining relation
 \eqref{def1} of $\R$-polynomials,
we obtain that
\begin{align}
\R_{e,\, v_{n,i}}(q)=&\R_{e,\, \Omega_n  s_{n-3} s_{n-4}\cdots s_{i-1}}(q)\nonumber\\[5pt]
=&\R_{s_{i-1},\, \Omega_n  s_{n-3} s_{n-4}\cdots  s_{i}}(q)+q\R_{e,\, \Omega_n  s_{n-3} s_{n-4}\cdots  s_{i}}(q)\nonumber\\[5pt]
=&\left(\R_{s_{i-1}s_i,\, \Omega_n  s_{n-3} s_{n-4}\cdots s_{i+1}}(q)+q\R_{s_{i-1},\, \Omega_n  s_{n-3} s_{n-4}\cdots s_{i+1}}(q)\right)\nonumber\\[5pt]
&\quad +q\left(\R_{s_{i},\, \Omega_n  s_{n-3} s_{n-4}\cdots  s_{i+1}}(q)
+q\R_{e,\, \Omega_n  s_{n-3} s_{n-4}\cdots  s_{i+1}}(q)\right)\nonumber\\[5pt]
=&\cdots\nonumber\\
=&\sum_{i-1\leq i_1<\cdots<i_k\leq n-3}q^{n-i-1-k}\R_{s_{i_1}\cdots  s_{i_k}, \, \Omega_n}(q).\label{0930-2}
\end{align}
Observe that $s_{i_1}\cdots  s_{i_k}$ is a reduced subword of $\Omega_n$
with $d(s_{i_1}\cdots  s_{i_k})=0$. By  Corollary \ref{coro-1},  we find that
\begin{equation}\label{ll}
\R_{s_{i_1}\cdots  s_{i_k}, \, \Omega_n}(q)=q^kT_{n-k}(q).
\end{equation}
Substituting \eqref{ll} into \eqref{0930-2}, we get
\begin{align*}
\R_{e,\, v_{n,i}}(q)=&\sum_{i-1\leq i_1<\cdots<i_k\leq n-3}q^{n-i-1-k}\R_{s_{i_1}\cdots  s_{i_k}, \, \Omega_n}(q)\\[5pt]
=& \sum_{i-1\leq i_1<\cdots<i_k\leq n-3}q^{n-i-1}\,T_{n-k}(q)\\[5pt]
=& \sum_{k=0}^{n-i-1}q^{n-i-1}{n-i-1 \choose k}T_{n-k}(q),
\end{align*}
as required.
\qed

We conclude this paper with  the following conjecture, which has been verified for $n\leq 9$.

\begin{conjecture}
For $n\geq 2$ and  $e\leq \sigma_1\leq \sigma_2\leq \Omega_n$, we have
\begin{equation}
\R_{\sigma_1,\sigma_2}(q)=
q^{g(\sigma_1,\sigma_2)}
\prod_{i=1}^{k}F_{h_i(\sigma_1,\sigma_2)}(q^{-2}),
\end{equation}
where $k$, $g(\sigma_1,\sigma_2)$ and $h_i(\sigma_1,\sigma_2)$
are integers depending on  $\sigma_1$
and $\sigma_2$.
\end{conjecture}

\vspace{.2cm} \noindent{\bf Acknowledgments.}
This work was
supported by the 973 Project  and the National Science Foundation of China.

\end{document}